\documentclass{svmult}
\usepackage{makeidx}         
\usepackage{graphicx}        
\usepackage{multicol}        
\usepackage[bottom]{footmisc}

\makeindex             

\newtheorem{thm}{Theorem}[section]
\newtheorem{lem}{Lemma}[section]
\newtheorem{prop}{Proposition}[section]
\newtheorem{cor}{Corollary}[section]

\newtheorem{defi}{Definition}[section]

\begin{document}

\title*{A Generalized It$\hat {\rm o}$'s Formula  in Two-Dimensions and Stochastic Lebesgue-Stieltjes Integrals}
\titlerunning{Two-dimensional generalized It$\hat {\rm o}$ Formula }
\author{Chunrong Feng\inst{1,2},
Huaizhong Zhao\inst{1}}
\authorrunning{C. Feng and H. Zhao}
\institute{ Department of Mathematical Sciences, Loughborough
University, LE11 3TU, UK. \texttt{C.Feng@lboro.ac.uk},
\texttt{H.Zhao@lboro.ac.uk}
\and School of Mathematics and System Sciences, Shandong University, Jinan, Shandong Province, 250100, China}

\maketitle
\newcounter{bean}
\begin{abstract}

In this paper, a generalized It${\hat {\rm o}}$ formula for time dependent functions of two-dimensional
    continuous semi-martingales is proved. The formula uses the local time of each
    coordinate process of the semi-martingale, left space and time first derivatives and second
    derivative
    $\nabla _1^- \nabla _2^-f$ only which are assumed to be of locally bounded variation
    in certain variables, and stochastic Lebesgue-Stieltjes integrals of two parameters.
    The two-parameter integral is defined as a natural generalization of the It${\hat {\rm o}}$ integral and Lebesgue-Stieltjes integral through a type of It${\hat {\rm o }}$ isometry formula.

 \vskip5pt

\end{abstract}

 \renewcommand{\theequation}{\arabic{section}.\arabic{equation}}

 \section{Introduction}

        The classical It$\hat {\rm o}$'s formula for twice differentiable functions
 has been extended to less smooth functions by many mathematicians.
 Progresses have been made mainly in one-dimension beginning with Tanaka's
 pioneering work \cite{tan} for $|X_t|$ to which the local time was beautifully linked.
  Further extensions were made to time independent convex
  functions in \cite{meyer} and \cite{wang}; to the case of absolutely
  continuous function with the first derivative being locally bounded
  in \cite{bou}; to $W_{loc}^{1,2}$ functions of a Brownian motion in \cite{Protter}  for one dimension and \cite{Protter2}
  for multi-dimensions.  It was proved in \cite{Protter} that $f(B_t)=f(B_0)+\int_0^t f'(B_s)dB_s+{1\over 2}[f(B),B]_t$,
  where $[f(B),B]_t$ is the covariation of the processes $f(B)$ and $B$ and is equal to $\int_0^t f(B_s)d^*B_s-\int_0^t f(B_s)dB_s$
  as a difference of backward and forward integrals. See \cite{rv} for the case of continuous semi-martingale. The multi-dimensional case was
  considered by \cite{Protter2}, \cite{rv} and \cite{mn}. An integral $\int_{-\infty}^\infty f^{\prime}(x){\rm d}_x L_t(x)$
  was introduced  in \cite {bou} through the existence of the expression $f(X(t))-f(X(0))-\int_0^t {\partial ^-\over \partial x} f(X(s))dX(s)$,
  where $L _t(x)$ is the local time of the semi-martingale $X_t$.  This work was extended further
   to define $\int_0^t\int_{-\infty}^\infty {\partial \over \partial x} f(s,X(s))d_{s,x}L_s(x)$ for a time dependent
   function $f(s,x)$ using forward and backward integrals for Brownian motion in \cite{eisenbaum1} and to semi-martingales
   other than Brownian motion in
   \cite{eisenbaum2}. This integral was also defined in \cite{rog}
   as a stochastic integral with excursion fields, and in \cite{Peskir1}
   through It$\hat {\rm o}$'s formula without assuming the reversibility of the semi-martingale which was required in \cite{eisenbaum1}.
   Other generalizations include \cite{frw} where it was also proved that if $X$ is a semi-martingale, then $f(X(t))$ is a semi-martingale if and only if $f\in W_{loc}^{1,2}$ and its weak derivative is of bounded variation using backward and forward integrals (\cite{lyons}).

The above mentioned extensions are useful in many problems. However, to use probabilistic methods to study problems arising in partial differential equations with singularities and mathematics of finance, we often need a generalized It$\hat {\rm o}$'s formula for time dependent $f(t,x)$.
In a special case that is if there exists a Radon measure $\nu$ and locally bounded Borel function $H$ such that $d_x(\nabla f(t,x))=H(t,x)\nu(dx)$, a generalized It$\hat {\rm o}$'s formula was obtained by \cite{yor1}. In a recent work \cite{Zhao1}, a new generalized It$\hat {\rm o}$'s formula for
 one-dimensional
continuous semi-martingales was proved. It is given in terms of a
Lebesgue-Stieltjes integral of the local time $L_t(x)$ with
respect to the two-dimensional variation of $\nabla^-f(t,x)$ as
follows
 \begin{eqnarray}\label{zhao1}
&&f(t,X(t))-f(0,X(0))\nonumber\\
&=&\int _0^t{\partial ^-\over \partial s} f(s,X(s)){\rm d}s+\int
_0^t\nabla ^-f(s,X(s))dX_s\nonumber\\
&&+{1\over 2} \int_0^t \Delta f_h(s,X(s))d<\hskip-4pt X\hskip-4pt>_s
+
\int _{-\infty}^{\infty }L _t(x){\rm d}_x\nabla ^-f_v(t,x)\nonumber\\
&&-\int
_{-\infty}^{+\infty}\int _0^{t}L _s(x)
{\bf \rm d}_{s,x}\nabla ^-f_v(s,x).\ \ a.s.
\end{eqnarray}
Here $f(t,x)=f_h(t,x)+f_v(t,x)$ is left continuous with $f_h(t,x)$ being $C^1$ in $x$ and
 $\nabla f_h(t,x)$ being absolutely continuous whose left derivative
 $\Delta ^-f_h(t,x)$ is left continuous and locally bounded, and
 $\nabla ^-f_v(t,x)$ being of locally bounded variation in $(t, x)$ and of locally bounded variation
 in $x$ at $t=0$.
Note the last two integrals are pathwise well defined due to the
well-known fact that the local time $L _t(x)$ is jointly
continuous in $t$ and c$\grave{a}$dl$\grave{a}$g in $x$ and has a
compact support in space $x$ for each $t$ (\cite{yor}, \cite{ks}).
In a special case, when there exists a curve $x=\gamma (t)$ of
locally bounded variation and the function $f$ is continuous but the
first order derivative $\nabla f$ has jumps across the curve and
second order derivative $\Delta f$ has left limit when $x\to
\gamma(t) -$, i.e. $\Delta ^-f$ exists and locally bounded and
left continuous off the curve(s) $x=\gamma (t)$, and there may be
jumps of $\nabla f$ along $x=\gamma (t)$ ($\Delta f$ is still
undefined), define $\Delta^-f$ on
the curve $x=\gamma(t)$
as the left limit of $\Delta f$. Then the following formula was derived
from (\ref{zhao1}) using the integration by parts formula
(\cite{Zhao1}):
\begin{eqnarray}\label{zhao2}
f(t,X(t))
&=&f(0,X(0))+\int _0^t{\partial ^-\over \partial s} f(s,X(s)){\rm d}s+\int
_0^t\nabla ^-f(s,X(s))dX_s\nonumber\\
&&
+
{1\over 2} \int _0^t\Delta ^-f(s,X(s))d<\hskip-4pt X\hskip-4pt>_s\nonumber\\
&& +\int _0^{t}(\nabla f(s,\gamma (s)+)-\nabla f(s,\gamma (s)-))dL
_s(\gamma(s)).\ \ a.s.
\end{eqnarray}
Here $dL_s(a)$ refers to the Lebesgue-Stieltjes integral with
respect to $s\mapsto L _s(a)$. Formula (\ref {zhao2}) was also
observed in \cite{peskir2} independently. These two new formulae
have been proved useful in analysing asymptotics of solutions of
partial differential equations in the presence of caustics
(\cite{Zhao2}) and studying the smooth fitting problem in American
put option (\cite{peskir4}). Formula (\ref{zhao1}) is in a very
general form. It includes the classical It${\hat{\rm o}}$ formula,
Tanaka's formula, Meyer's formula for convex functions, the
formula given by Az\'ema, Jeulin, Knight and Yor \cite{yor} and
formula (1.2).

The purpose of this paper is to extend formula (\ref{zhao1}) to two dimensions. This is
a nontrivial extension as the  local time in two-dimensions does not exist.
But we observe
for a smooth function $f$, formally by the occupation times formula
\begin{eqnarray}\label{zhao3}
&& {1\over 2}\int _0^{t}\Delta _1f(s,X_1(s),X_2(s))d<\hskip-4pt X
_1\hskip-4pt>_s
\nonumber\\
&=& \int
_{-\infty}^{+\infty}\int _0^{t}\Delta _1f(s,a,X_2(s)){\rm d}_sL
_1(s,a){\rm d}a\nonumber\\
&=&
\int _{-\infty}^{+\infty}\Delta _1f({t},a,X_2({t})){L_1(t,a)}{\rm d}a\nonumber\\
&&-\int
_{-\infty}^{+\infty}\int _0^{t}{L_1(s,a)}{\bf\rm d}_{s,a}\nabla_1 f(s,a,X_2(s)),
\end{eqnarray}
if the integral $\int_{-\infty}^{+\infty}\int
_0^{t}{L_1(s,a)}{\bf\rm d}_{s,a}\nabla _1f(s,a,X_2(s))$ is
properly defined. Here $\nabla_1 f(s,a,X_2(s))$ is a
semi-martingale for any fixed $a$, following the one-dimensional
generalized It${\hat {\rm o}}$'s formula (\ref{zhao1}). For this,
we study this kind of the integral $\int_{-\infty}^{+\infty}\int
_0^{t}{g(s,a)}{\bf\rm d}_{s,a} h(s,a)$ in section 2.  Here
$h(s,x)$ is a continuous martingale with cross variation
$<h(\cdot,a),h(\cdot,b)>_s$ of locally bounded variation in
$(a,b)$, and $E\left[\int_0^t\int_{R^2}|g(s,a)g(s,b)||{\bf \rm
d}_{a,b,s}<h(\cdot,a),h(\cdot,b)>_s|\right] < {\infty}$. The
integral is different from the Lebesgue-Stieltjes integral and
It${\hat{\rm o}}$'s stochastic integral. But it is a natural
extension to the two-parameter stochastic case and therefore
called a stochastic Lebesgue-Stieltjes integral. According to our
knowledge, our integral is new. It's different from integration
with Brownian sheet defined by Walsh (\cite{walsh}) and
integration w.r.t. Poisson random measure (see \cite{ikeda}). A
generalized It$\hat {\rm o}$'s formula in two dimensions is proved
in section 3. Applications e.g. in the study of the asymptotics of
the solutions of heat equations with caustics in two dimensions,
are not included in this paper. These results will be published in
some future work. Furthermore, it has been observed by us in
\cite{Zhao3}  that the local time $L_t(x)$ can be considered
naturally as a rough path in $x$ of finite 2-variation and
$\int_0^t\int_{-\infty}^\infty \nabla^-f(s,x)d_{s,x}L_s(x)$ is
defined pathwisely by using and extending Lyons' idea of rough
path integration (\cite{terry}). \vskip5pt

\section{The definition of stochastic Lebesgue-Stieltjes integrals and the integration by parts formula}
\setcounter{equation}{0}

For a filtered probability space $(\Omega,{\cal F},\{{\cal F}_t\}_{t\geq 0},P)$,
denote by ${\cal M}_2$ the Hilbert space of all processes $X=(X_t)_{0\leq t\leq T}$ such that
 $(X_t)_{0\leq t\leq T}$ is a $({\cal F}_t)_{0\leq t\leq T}$ right continuous square integrable martingale with
 inner product $(X,Y)=E(X_TY_T)$. A three-variable function $f(s,x,y)$ is called left continuous iff it is
 left continuous in all three variables together i.e. for any sequence $(s_1,x_1,y_1)\leq   (s_2,x_2,y_2)
 \leq \cdots \leq (s_k,x_k,y_k)\to (s,x,y)$, as $k\to \infty$, we have $f(s_k,x_k,y_k)\to f(s,x,y)$ as $k\to \infty$. Here $(s_1,x_1,y_1)\leq
 (s_2,x_2,y_2)$ means $s_1\leq s_2$, $x_1\leq x_2$ and $y_1\leq y_2$.
 Define
\begin{eqnarray*}
{\cal V}_1:= \Big\{h: && [0,t] \times ({-\infty},{\infty}) \times {\Omega} \to R \ s.t.\ (s,x,\omega)\mapsto h(s,x,\omega)\\
&& is \  {\cal B} ([0,t] \times R) \times {\cal F}{\rm-}measurable,\ and\ h(s,x)\  is\  \\
&&{\cal F}_s{\rm-}adapted \ for\ any\ x\in R\Big \},\\
 {\cal V}_2:=\Big\{h: && h\in{\cal V}_1\  is\  a \ continuous \ (in\  s)\ {\cal M}_2-martingale
 \  for \ each \ x,\\
 && and \ the \ crossvariation\ <h(\cdot,x),h(\cdot,y)>_s is\ left\ continuous\\
 && and \ of \ locally\ bounded\  variation \ in\  (s,x,y)\Big\}.
\end{eqnarray*}
In the following, we will always denote $<h(\cdot,x),h(\cdot,y)>_s$ by $<h(x),h(y)>_s$.

We now recall some classical results for the sake of completeness
of the paper (see \cite{ash} and \cite{mc}). A three-variable
function $f(s,x,y)$ is called monotonically increasing if whenever
$(s_2,x_2,y_2)\geq (s_1,x_1,y_1)$, then
\begin{eqnarray*}
&&f(s_2,x_2,y_2)-f(s_2,x_1,y_2)-f(s_2,x_2,y_1)+f(s_2,x_1,y_1)\\
&&-
f(s_1,x_2,y_2)+f(s_1,x_1,y_2)+f(s_1,x_2,y_1)-f(s_1,x_1,y_1)\geq 0.
\end{eqnarray*}
For a left-continuous and monotonically increasing function $f(s,x,y)$, one can define a
Lebesgue-Stieltjes measure by setting
\begin{eqnarray*}
&&
\nu ([s_1,s_2)\times [x_1,x_2)\times [y_1,y_2))\\
&=&f(s_2,x_2,y_2)-f(s_2,x_1,y_2)-f(s_2,x_2,y_1)+f(s_2,x_1,y_1)\\
&&-
f(s_1,x_2,y_2)+f(s_1,x_1,y_2)+f(s_1,x_2,y_1)-f(s_1,x_1,y_1).
\end{eqnarray*}
For $h\in{\cal V}_2$, define
\begin{eqnarray*}
<h(x),h(y)>_{t_1}^{t_2}:=<h(x),h(y)>_{t_2}
- <h(x),h(y)>_{t_1},\ {t_2}\geq{t_1}.
\end{eqnarray*}
 Note as $<h(x),h(y)>_s$ is left continuous and of locally bounded variation in ($s,x,y$), so it can be
 decomposed to the difference of two increasing and left continuous functions $f_1(s,x,y)$ and $f_2(s,x,y)$
(see McShane \cite{mc} or Proposition 2.2 in Elworthy, Truman and
Zhao \cite{Zhao1} which also holds for multi-parameter functions).
Note each of $f_1$ and $f_2$ generates a measure, so for any
measurable function $g(s,x,y)$, we can define
\begin{eqnarray*}
&&
\int _{t_1}^{t_2}\int _{a_1}^{a_2}\int _{b_1}^{b_2}g(s,x,y)d_{x,y,s}<h(x),h(y)>_s\\
&=& \int _{t_1}^{t_2}\int _{a_1}^{a_2}\int
_{b_1}^{b_2}g(x,y,s)d_{x,y,s}f_1(s,x,y)\\
&&- \int _{t_1}^{t_2}\int _{a_1}^{a_2}\int
_{b_1}^{b_2}g(x,y,s)d_{x,y,s}f_2(s,x,y).
\end{eqnarray*}
In particular, a signed product measure in the space $[0,T]\times R^2$ can be defined as follows: for any ${[t_1,t_2)}\times{[x_1,x_2)}\times{[y_1,y_2)} \subset [0,T]\times R^2$
\begin{eqnarray}
&&
 \int _{t_1}^{t_2}\int _{x_1}^{x_2} \int _{y_1}^{y_2} {\bf \rm d}_{x,y,s}<h(x),h(y)>_s\nonumber\\
&=&\int _{t_1}^{t_2}\int _{x_1}^{x_2} \int _{y_1}^{y_2}{\bf \rm
d}_{x,y,s}f_1(s,x,y)-\int _{t_1}^{t_2}\int _{x_1}^{x_2} \int
_{y_1}^{y_2}{\bf \rm d}_{x,y,s}f_2(s,x,y)\nonumber\\
 &=&<h(x_2),h(y_2)>_{t_1}^{t_2} - <h(x_2),h(y_1)>_{t_1}^{t_2}\nonumber\\
 &&-<h(x_1),h(y_2)>_{t_1}^{t_2} + <h(x_1),h(y_1)>_{t_1}^{t_2}\nonumber\\
 &=&<h(x_2)-h(x_1),h(y_2)-h(y_1)>_{t_1}^{t_2}.
 \end{eqnarray}
 Define
 \begin{eqnarray}
 |{\bf \rm d}_{x,y,s}<h(x),h(y)>_s|={\bf \rm d}_{x,y,s}f_1(s,x,y)+{\bf \rm d}_{x,y,s}f_2(s,x,y).
 \end{eqnarray}
 Moreover, for $h\in{\cal V}_2$, define:
 \begin{eqnarray*}
 {\cal V}_3(h):=\Big\{g&:&g\in{\cal V}_1\  has\  a\  compact\  support\  in\  x,\  and\
   \\&E&\left[\int_0^t\int_{R^2}|g(s,x)g(s,y)||{\bf \rm d}_{x,y,s}<h(x),h(y)>_s|\right] < {\infty} \Big\}.
 \end{eqnarray*}

Consider now a simple function
\begin{eqnarray}\label{fr1}
  g(s,x,\omega)&=&\sum_{i=1}^{n}\sum_{j=1}^{m}e(t_j,x_i)1_{(t_j,t_{j+1}]}(s)1_{(x_i,x_{i+1}]}(x)
\end{eqnarray}
  where $t_1<t_2<\cdots <t_{n+1}$, $x_1<x_2<\cdots <x_{m+1}$, $ e(t_j,x_i)$ are ${\cal F}_{t_j}$-measurable. For $h\in{\cal  V}_2$, define an integral as:
\begin{eqnarray}
I_t(g)&:=&\int _{0}^{t}\int _{-\infty}^{\infty}g(s,x){\bf \rm d}_{s,x}h(s,x)\nonumber\\
&=&\sum _{i=1}^{n}\sum _{j=1}^{m} e(t_j\wedge t,x_i)\Big[h(t_{j+1}\wedge t,x_{i+1})-h(t_j\wedge t,x_{i+1})\nonumber\\
&&\hskip 3cm -h(t_{j+1}\wedge t,x_i)+h(t_j\wedge t,x_i)\Big].
\end{eqnarray}
This integral is called the  stochastic Lebesgue-Stieltjes
integral of the simple function $g$. It's easy to see for simple
functions $g_1, g_2\in {\cal V}_3(h)$,
\begin{eqnarray}\label{add}
I_t(\alpha g_1+\beta g_2)=\alpha I_t(g_1)+\beta I_t(g_2),
\end{eqnarray}
for any $\alpha, \beta\in R$.
 The following lemma plays a
key role in extending the integral of simple functions to
functions in ${\cal V}_3(h)$. It is equivalent to the It$\hat {\rm
o}$'s isometry formula in the case of the stochastic integral.
\vskip5pt

\begin{lem}\label{lemma1}
 If $h\in{\cal V}_2$, $g\in{\cal V}_3(h)$ is simple, then $I_t(g)$ is a continuous martingale
 with respect to $({\cal F}_t)_{0\leq t\leq T}$ and
\begin{eqnarray}\label{fcr1}
 &&E \left(\int _{0}^{t}\int _{-\infty}^{\infty}g(s,x){\bf \rm d}_{s,x}h(s,x)\right)^2\nonumber\\
 &=&E\int _{0}^{t}\int _{R^2} g(s,x)g(s,y){\bf \rm d}_{x,y,s}<h(x),h(y)>_s.
\end{eqnarray}
\end{lem}
 {\em Proof}: From the definition of $\int _{0}^{t}\int _{-\infty}^{\infty}g(s,x){\bf \rm d}_{s,x}h(s,x)$,
it is easy to see that $I_t$ is a continuous martingale
 with respect to $({\cal F}_t)_{0\leq t\leq T}$.
 As $h(s,x,\omega)$ is a continuous martingale in ${\cal M}_2$,
 using a standard conditional expectation argument
 to remove the cross product parts, we get:
\begin{eqnarray*}
&&
E\left[\left(\int _{0}^{t}\int _{-\infty}^{\infty}g(s,x){\bf \rm d}_{s,x}h(s,x)\right)^2\right]\\
&=&E\sum _{j=1}^{m}\Bigg(\sum _{i=1}^{n} e(t_j\wedge t,x_i)\Big[h(t_{j+1}\wedge t,x_{i+1})-h(t_{j}\wedge t,x_{i+1})\nonumber\\
&&\hskip 3cm -h(t_{j+1}\wedge t,x_i)+h(t_{j}\wedge t,x_{i})\Big]\Bigg)^2 \\
&=&E\sum _{j=1}^{m}\Bigg (\sum _{i=1}^{n} \sum _{k=1}^{n} e(t_j\wedge t,x_i)e(t_j\wedge t,x_k)\cdot\\
&&\hskip 1cm \Big[h(t_{j+1}\wedge t,x_{i+1}) - h(t_{j}\wedge t,x_{i+1}) - h(t_{j+1}\wedge t,x_i) + h(t_{j}\wedge t,x_{i})\Big]\cdot \\
&&\hskip 1cm\Big[h(t_{j+1}\wedge t,x_{k+1}) - h(t_j\wedge t,x_{k+1}) - h(t_{j+1}\wedge t,x_k) + h(t_j\wedge t,x_k)\Big]\Bigg)\\
&=&E\sum _{j=1}^{m}\Bigg\{\sum _{i=1}^{n} \sum _{k=1}^{n}
e(t_j\wedge t,x_i)e(t_j\wedge t,x_k)\cdot\\
&&\hskip 1cm\Big[\big(h(t_{j+1}\wedge t,x_{i+1})-h(t_{j}\wedge t,x_{i+1})\big)\big(h(t_{j+1}\wedge t,x_{k+1})-h(t_{j}\wedge t,x_{k+1})\big)\\
&& \hskip 1cm- \big(h(t_{j+1}\wedge t,x_{i+1})-h(t_{j}\wedge t,x_{i+1})\big)\big(h(t_{j+1}\wedge t,x_k) - h(t_{j}\wedge t,x_k)\big)  \\
&& \hskip 1cm - \big(h(t_{j+1}\wedge t,x_i) - h(t_{j}\wedge t,x_{i})\big)\big(h(t_{j+1}\wedge t,x_{k+1})-h(t_{j}\wedge t,x_{k+1})\big) \\
&&\hskip 1cm + \big(h(t_{j+1}\wedge t,x_i) - h(t_{j}\wedge t,x_{i})\big)\big(h(t_{j+1}\wedge t,x_k) - h(t_{j}\wedge t,x_k)\big)\Big]\Bigg\}\\
&=&E \int _{0}^{t} \sum _{i=1}^{n} \sum _{k=1}^{n}e(s,x_i)e(s,x_k)
\Big[{\bf \rm d}_s <h(x_{i+1}),h(x_{k+1})>_s \\
&&- {\bf \rm d}_s <h(x_{i+1}),h(x_{k})>_s
- {\bf \rm d}_s <h(x_{i}),h(x_{k+1})>_s \\
&&
\hskip3.8cm + {\bf \rm d}_s <h(x_{i}),h(x_{k})>_s\Big]\\
&=&E\int _{0}^{t} \sum _{i=1}^{n} \sum _{k=1}^{n}e(s,x_i)e(s,x_k)) \Big[{\bf \rm d}_s <h(x_{i+1})-h(x_{i}),h(x_{k+1})-h(x_{k})>_s\Big]\\
&=&E\left[\int _{0}^{t}\int _{R^2} g(s,x)g(s,y){\bf \rm
d}_{x,y,s}<h(x),h(y)>_s\right].
\end{eqnarray*}
So we prove the desired result.
$\hfill\diamond$\\

The idea is to use (\ref {fcr1}) to extend the definition of the
integrals of simple functions to integrals of functions in ${\cal
V}_3(h)$, for any $h \in {\cal V}_2$. We achieve this goal in
several steps:

\begin{lem}\label{lem1}   Let $h \in {\cal V}_2$, $f \in {\cal V}_3(h)$ be bounded uniformly in $\omega$,
$f(\cdot,\cdot,\omega)$ be continuous for each ${\omega}$ on its
compact support. Then there exist a sequence of bounded simple
functions  ${\varphi}_{m,n} \in {\cal V}_3(h) $  such that
\begin{eqnarray*}
&&\hskip -0.3cm E\int _{0}^{t}\int _{R^2} \mid (f - \varphi _{m,n})(s,x)(f - \varphi_{m',n'})(s,y) |\mid{\bf \rm d}_{x,y,s}<h(x),h(y)>_s|\\
&&\hskip -0.3cm \to0,
\end{eqnarray*}
as $m, n , m',n'\to \infty$.
\end{lem}
{\em Proof}: Let $[0,T]\times[a,b]$ be a rectangle covering the
compact support of $f$ and $0=t_1<t_2<\cdots<t_{m+1}=T$, and
$a=x_1<x_2<\cdots<x_{n+1}=b$ be a partition of $[0,T]\times[a,b]$.
Assume when $n,m\to\infty$, $\max \limits _{1\leq j\leq
m}(t_{j+1}-t_j)\to 0$, $\max \limits _{1\leq i\leq
n}(x_{i+1}-x_i)\to 0$. Define
\begin{eqnarray}
\varphi _{m,n}(t,x):=\sum_{j=1}^{m}\sum_{i=1}^{n}f(t_j,x_i)1_{(t_j,t_{j+1}]}(t)1_{(x_i,x_{i+1}]}(x).
\end{eqnarray}
Then $\varphi _{m,n}(t,x)$ are simple and $\varphi _{m,n}(t,x)\to f(t,x) \ a.s. \ as \ m,n\to {\infty}$.
The result follows Lebesgue's dominated convergence theorem.
$\hfill\diamond$

\begin{lem} \label{lem2}  Let $h \in {\cal V}_2$ and $k \in {\cal V}_3(h)$ be bounded uniformly in $\omega$.
Then there exist  functions $f_n \in {\cal V}_3(h) $ such that $f_n(\cdot,\cdot,\omega) $ are continuous
for all $\omega$ and $n$ on its support, and
\begin{eqnarray*}
&&E \int _{0}^{t}\int _{R^2} \mid (k - f_n)(s,x)(k - f_{n'})(s,y)
|\mid{\bf \rm d}_{x,y,s}<h(x),
h(y)>_s|\\
&&\to0,
\end{eqnarray*}
as $n,n' \to {\infty}$.
\end{lem}
{\em Proof}: Define
\begin{eqnarray*}
f_{n}(s,x)=n^2\int_{x-{1\over n}}^x\int_{s-{1\over n}}^s k(\tau,y)d\tau dy.
\end{eqnarray*}
Then $f_{n}(s,x)$ is continuous in $s, x$, and when $n \to {\infty} $, $f_{n}(s,x)\to k(s,x)$ a.s.. The desired convergence follows by Lebesgue's dominated convergence theorem.
$\hfill\diamond$

\begin{lem} \label{lem3}  Let $h \in {\cal V}_2$ and $g\in{\cal V}_3(h)$. Then there exist  functions $k_n\in {\cal V}_3(h)$, bounded uniformly in $\omega$  for each $n$, and
\begin{eqnarray*}
&&E\int _{0}^{t}\int _{R^2} \mid (g - k_n)(s,x)(g - k_{n'})(s,y)
|\mid{\bf \rm d}_{x,y,s}<h(x),
h(y)>_s|\\
&&\to0,
\end{eqnarray*}
as $ n,n' \to {\infty}$.
\end{lem}
{\em Proof}: Define
\begin{eqnarray}
 k_n(t,x,\omega):=\cases {-n   { \ \ \ \ \ \ \ \  \ \ {\rm  if}\ g(t,x,\omega)<-n} \cr g(t,x,\omega)   { \ \  {\rm  if} \ -n \leq g(t,x,\omega) \leq n} \cr n   { \ \ \ \ \ \ \ \ \ \ \  \ {\rm  if} \ g(t,x,\omega)>n.}}
\end{eqnarray}
Then as $n\to \infty$, $k_n(t,x,\omega)\to g(t,x,\omega)$ for each
($t,x,\omega$). Note $|k_n(t,x,\omega)|\leq |g(t,x,\omega)|$ and
$g\in{\cal V}_3(h)$. So applying Lebesgue's dominated convergence
theorem, we obtain the desired result. $\hfill\diamond$ \vskip5pt

From Lemmas \ref{lem3}, \ref{lem2}, \ref{lem1}, for each $h \in
{\cal V}_2$, $g \in {\cal V}_3(h)$, we can construct a sequence of
simple functions $\{{\varphi} _{m,n}\} $ in ${\cal V}_3(h)$ such
that,
\begin{eqnarray*}
&&E\int _{0}^{t}\int _{R^2} \mid {(g - \varphi _{m,n})(s,x)(g - \varphi _{m',n'})(s,y)} |\mid{\bf \rm d}_{x,y,s}<h(x),h(y)>_s|\\
&&\to0,
\end{eqnarray*}
as $ m,n,m',n' \to {\infty}$. For $\varphi _{m,n}$ and $\varphi
_{m',n'}$, we can define stochastic Lebesgue-Stieltjes integrals
$I_t(\varphi _{m,n})$ and $I_t(\varphi _{m',n'})$. From Lemma
\ref{lemma1} and (\ref{add}), it is easy to see that
\begin{eqnarray*}
&&E\left[I_T(\varphi _{m,n})-I_T(\varphi
_{m',n'})\right]^2\\
&=&E\left[I_T(\varphi _{m,n}-\varphi
_{m',n'})\right]^2\\
&=&E\int_0^T\int_{R^2}(\varphi
_{m,n}-\varphi_{m',n'})(s,x)(\varphi
_{m,n}-\varphi_{m',n'})(s,y)d_{x,y,s}<h(x),h(y)>_s\\
&=&E\int_0^T\int_{R^2}[(\varphi
_{m,n}-g)-(\varphi_{m',n'}-g)](s,x)\cdot\\
&&\hskip 1.5cm[(\varphi
_{m,n}-g)-(\varphi_{m',n'}-g)](s,y)d_{x,y,s}<h(x),h(y)>_s\\
&=&E\int_0^T\int_{R^2}(\varphi _{m,n}-g)(s,x)(\varphi
_{m,n}-g)(s,y)d_{x,y,s}<h(x),h(y)>_s \\
&&-E\int_0^T\int_{R^2}(\varphi _{m,n}-g)(s,x)(\varphi
_{m',n'}-g)(s,y)d_{x,y,s}<h(x),h(y)>_s \\
&&-E\int_0^T\int_{R^2}(\varphi _{m',n'}-g)(s,x)(\varphi
_{m,n}-g)(s,y)d_{x,y,s}<h(x),h(y)>_s \\
&&+E\int_0^T\int_{R^2}(\varphi _{m',n'}-g)(s,x)(\varphi
_{m',n'}-g)(s,y)d_{x,y,s}<h(x),h(y)>_s \\
&\leq&E\int_0^T\int_{R^2}\mid(\varphi _{m,n}-g)(s,x)(\varphi
_{m,n}-g)(s,y)\mid\mid d_{x,y,s}<h(x),h(y)>_s\mid \\
&&+E\int_0^T\int_{R^2}\mid(\varphi _{m,n}-g)(s,x)(\varphi
_{m',n'}-g)(s,y)\mid\mid d_{x,y,s}<h(x),h(y)>_s\mid \\
&&+E\int_0^T\int_{R^2}\mid(\varphi _{m',n'}-g)(s,x)(\varphi
_{m,n}-g)(s,y)\mid\mid d_{x,y,s}<h(x),h(y)>_s\mid \\
&&+E\int_0^T\int_{R^2}\mid(\varphi _{m',n'}-g)(s,x)(\varphi
_{m',n'}-g)(s,y)\mid\mid d_{x,y,s}<h(x),h(y)>_s \mid\\
&\to& 0,
\end{eqnarray*}
as $m,n,m',n'\to \infty$. Therefore
$\{I_.(\varphi_{m,n})\}_{m,n=1}^\infty$ is a Cauchy sequence in
${\cal M}_2$ whose norm is denoted by $\parallel\cdot\parallel$.
So there exists a process $I(g)=\{I_t(g), 0\leq t\leq T\}$ in
${\cal M}_2$, defined modulo indistinguishability, such that
\begin{eqnarray*}
\parallel I(\varphi_{m,n})-I(g)\parallel\to 0, \ as \
m,n\to\infty.
\end{eqnarray*}
By the same argument as for the stochastic integral, one can
easily prove that  $I(g)$ is well-defined (independent of the choice of the simple functions), and (\ref{fcr1}) is true for  $I(g)$. We now can have the following
definition.

\begin{defi}\label{definition1}
Let $h \in {\cal V}_2$, $g \in {\cal V}_3(h)$.Then the integral of
$g$ with respect to $h$ can be defined as:
\begin{eqnarray*}
&&\int _{0}^{t}\int _{-\infty}^{\infty}g(s,x){\bf \rm d}_{s,x}h(s,x)\\
&=&\lim_{m,n \to {\infty}} \int _{0}^{t}\int _{-\infty}^{\infty}\varphi_{m,n}(s,x){\bf \rm d}_{s,x}h(s,x), \ \ \ \ (limit \ in \ {\cal M}_2)
\end{eqnarray*}
is a continuous martingale with respect to $({\cal F}_t)_{0\leq t\leq T}$ and for each $t\leq T$,
(\ref{fcr1}) is satisfied.
Here $\{{\varphi} _{m,n}\} $ is a sequence of simple functions in ${\cal V}_3(h)$, s.t.
\begin{eqnarray*}
&&E\int _{0}^{t}\int _{R^2} \mid (g - \varphi _{m,n})(s,x)(g - \varphi _{m',n'})(s,y) |\mid{\bf \rm d}_{x,y,s}<h(x),h(y)>_s|\\
&&\to0,
\end{eqnarray*}
as $ m,n, m^{\prime}, n^{\prime} \to {\infty}$. Note $\varphi _{m,n}$ may be constructed by combining the three approximation procedures in Lemmas \ref{lem3}, \ref{lem2}, \ref{lem1}.
\end{defi}

The following integration by parts formula will be useful in the
proof of our main theorem. Although the conditions are strong and
may be unnecessary, the proposition is enough for our purpose. We
don't strike to weaken the conditions here. \vskip5pt

\begin{prop}\label{proposition1}
 If $h \in {\cal V}_2$, $g \in {\cal V}_3(h)$, and $g(t,x)$ is $C^2$ in $x$, $\Delta g(t,x)$ is bounded uniformly in $t$, then a.s.
\begin{eqnarray}\label{fc1}
-\int _{-\infty}^{+\infty}\int _0^t \nabla g(s,x) {\rm d}_{s} h(s,x) dx=\int _0^t \int _{-\infty}^{+\infty}g(s,x){\rm \bf d}_{s,x} h(s,x).
\end{eqnarray}
\end{prop}
{\em Proof}:
If $g$ is a simple function as given in (\ref{fr1}),
one can always add some points in the partition to make $e(t_j\wedge t,x_1)=0$ and $e(t_j\wedge t,x_{n+1})=0$ for all $j=1,2,\cdots,m$ as $g$ has a compact support in $x$. So for $h\in{\cal V}_2$,
\begin{eqnarray*}
&&
\int _{0}^{t}\int _{-\infty}^{\infty}g(s,x){\bf \rm d}_{s,x}h(s,x)\\
&=&\sum _{i=1}^{n}\sum _{j=1}^{m} e(t_j\wedge t,x_i)\Big[h(t_{j+1}\wedge t,x_{i+1})-h(t_j\wedge t,x_{i+1})\\
&&\hskip 3cm -h(t_{j+1}\wedge t,x_i)+h(t_j\wedge t,x_i)\Big]\\
&=&-\sum _{i=0}^{n-1}\sum _{j=1}^{m} e(t_j\wedge t,x_{i+1})\Big[h(t_{j+1}\wedge t,x_{i+1})-h(t_j\wedge t,x_{i+1})\Big]\\
&&+\sum _{i=1}^{n}\sum _{j=1}^{m} e(t_j\wedge t,x_{i})\Big[h(t_{j+1}\wedge t,x_{i+1})-h(t_j\wedge t,x_{i+1})\Big]\\
&=&-\sum _{i=1}^{n}\sum _{j=1}^{m} \Big[e(t_j\wedge t,x_{i+1 })-e(t_j\wedge t,x_i)\Big]\Big[h(t_{j+1}\wedge t,x_{i+1})-h(t_j\wedge t,x_{i+1})\Big].
\end{eqnarray*}
If $g(t,x)$ is $C^2$ in $x$ and $\Delta g(t,x_2)$ is bounded uniformly in $t$,
let
\begin{eqnarray*}
\varphi _{m,n}(t,x):=\sum_{j=1}^{m}\sum_{i=1}^{n}g(t_j,x_i)1_{[t_j,t_{j+1})}(t)1_{[x_i,x_{i+1})}(x),
  \end{eqnarray*}
so
\begin{eqnarray*}
\varphi _{m,n}(t,x)\to {g(t,x)}\  a.s. \ as \ m,n\to {\infty}.
 \end{eqnarray*}
Then by the intermediate value theorem, there exist $\xi_i\in [x_i,x_{i+1}]$ $(i=1,2,\cdots,n)$ such that,
\begin{eqnarray*}
&&\int _{-\infty}^{+\infty}\int _0^t g(s,x) {\rm d}_{s,x}{h(s,x)}\\
&=&-\lim_{\delta_t,\delta_x \to 0}\sum _{i=1}^{n}\sum _{j=1}^{m} \Big[g(t_j\wedge t,x_{i+1})-g(t_j\wedge t,x_i)\Big]\\
&&\hskip 2.5cm\Big[h(t_{j+1}\wedge t,x_{i+1})-h(t_j\wedge t,x_{i+1})\Big] \ \ \ \ \ \ \ (limit \  in \ {\cal M}_2)\\
&=&-\lim_{\delta_t,\delta_x \to 0}\sum _{i=1}^{n}\sum _{j=1}^{m} \nabla g(t_j\wedge t,\xi_i)\Big[h(t_{j+1}\wedge t,x_{i+1})-h(t_j\wedge t,x_{i+1})\Big]\cdot\\
&&\hskip 2.5cm (x_{i+1}-x_i)\\
&=&-\lim_{\delta_x \to 0}\sum_{i=1}^{n}\int_0^t \nabla g(s,\xi_i){\rm d}_{s}h(s,x_{i+1})(x_{i+1}-x_{i})\ \ \ \ \ \ \ \ \ \ \ \ \  (limit \  in \ {\cal M}_2)\\
&=&-\lim_{\delta_x \to 0}\sum_{i=1}^{n}\int_0^t \nabla g(s,x_{i+1}){\rm d}_{s}h(s,x_{i+1})(x_{i+1}-x_{i})\\
&&-\lim_{\delta_x \to 0}\sum_{i=1}^{n}\int_0^t (\nabla g(s,\xi_i)-\nabla g(s,x_{i+1})){\rm d}_{s}h(s,x_{i+1})(x_{i+1}-x_{i})\\
&=&-\int _{-\infty}^{+\infty}\int _0^t \nabla g(s,x) {\rm d}_{s} h(s,x) dx.\hskip 3.5cm  (limit \  in \ {\cal M}_2)
\end{eqnarray*}
Here $ {\delta _t}=\max\limits_{1\leq j\leq m}{|t_{j+1}-t_j|} $, $  {\delta _x}=\max\limits_{1\leq i\leq m}{|x_{i+1}-x_i|} $. To prove the last equality, first notice that
\begin{eqnarray*}
&&\lim_{\delta_x \to 0}\sum_{i=1}^{n}\int_0^t \nabla g(s,x_{i+1}){\rm d}_{s}h(s,x_{i+1})(x_{i+1}-x_{i})\\
&=&\int _{-\infty}^{+\infty}\int _0^t \nabla g(s,x) {\rm d}_{s} h(s,x) dx.
\end{eqnarray*}
Second, by the intermediate value theorem again, the second term can be estimated as:
\begin{eqnarray*}
&&E\left[\sum_{i=1}^{n}\int_0^t (\nabla g(s,\xi_i)-\nabla g(s,x_{i+1})){\rm d}_{s}h(s,x_{i+1})(x_{i+1}-x_{i})\right]^2\\
&=&E \sum_{i=1}^{n}\sum_{k=1}^{n}\bigg[\int_0^t (\nabla g(s,\xi_i)-\nabla g(s,x_{i+1})){\rm d}_{s}h(s,x_{i+1})(x_{i+1}-x_{i})\cdot\\
&&\hskip 2cm \int_0^t (\nabla g(s,\xi_k)-\nabla g(s,x_{k+1})){\rm d}_{s}h(s,x_{k+1})(x_{k+1}-x_{k})\bigg]\\
&=&\sum_{i=1}^{n}\sum_{k=1}^{n}E\int_0^t (\nabla g(s,\xi_i)-\nabla g(s,x_{i+1}))(\nabla g(s,\xi_k)-\nabla g(s,x_{k+1}))\\
&&\hskip 2cm{\rm d}_s<h(x_{i+1}),h(x_{k+1})>_s(x_{i+1}-x_i)(x_{k+1}-x_{k})\\
&&\\
&\leq&\sum_{i=1}^{n}\sum_{k=1}^{n}E\sup\limits_{\xi_i\in[x_i,x_{i+1}]}|\nabla g(s,\xi_i)-\nabla g(s,x_{i+1})|\cdot\\
&&\hskip 1.5cm\sup\limits_{\xi_k\in[x_k,x_{k+1}]}|\nabla g(s,\xi_k)-\nabla g(s,x_{k+1})|\cdot\\
&&\hskip 1.5cm|<h(x_{i+1})>_t<h(x_{k+1})>_t|^{1\over2}(x_{i+1}-x_i)(x_{k+1}-x_k)\\
&\leq&E\Big [\sup\limits_{s}\sup\limits_{i}\sup\limits_{\xi_i\in[x_i,x_{i+1}]}|\Delta g(s,\eta_i)(\xi _i-x_{i+1})|\cdot\\
&&\hskip 0.3cm\sup\limits_{s}\sup\limits_{k}\sup\limits_{\xi_k\in[x_k,x_{k+1}]}|\Delta g(s,\eta_k)(\xi _k-x_{k+1})|
|<h(x_{i+1})>_t<h(x_{k+1})>_t|^{1\over2}\Big]\\
&&\hskip 0.3cm\cdot\left(\sum_{i=1}^{n}\sum_{k=1}^{n}(x_{i+1}-x_i)(x_{k+1}-x_k)\right)\\
&\to& 0,\ as\  \delta _x\to 0,
\end{eqnarray*}
    where $\eta _i\in [\xi _i, x_{i+1}]$, $\eta _k\in [\xi _k, x_{k+1}]$. The desired result is proved.
    \hfill $\diamond$
 \vskip5pt
\section{The generalized It${\hat{\rm o}}$'s formula in two-dimensional space}

\setcounter{equation}{0}

Let $X(s)=(X_1(s),X_2(s))$ be a two-dimensional continuous
semi-martingale with $X_i(s)=X_i(0)+M_i(s)+V_i(s) (i=1,2)$ on a
probability space $(\Omega,{\cal F},P)$. Here $M_i(s)$ is a
continuous local martingale and $V_i(s)$ is an adapted continuous
process of locally bounded variation (in $s$). Let $L_i(t,a)$ be
the local time of $X_i(t)$ (i=1,2)
\begin{eqnarray}
L _i(t,a)=\lim_{\epsilon\downarrow 0}
{1\over 2\epsilon}\int _0^t1_{[a,a+\epsilon)}(X_i(s))d<\hskip-4pt M_i\hskip-4pt>_s, \ \
a.s. \ \  i=1,2
\end{eqnarray}
for each $t$ and $a\in R$.
 Then it is well known for each
fixed $a\in R$,
  $L_i(t,a,\omega)$
is continuous, and nondecreasing in $t$ and right continuous with
left limit (c$\grave{a}$dl$\grave{a}$g) with respect to $a$
(\cite{ks}, \cite{yor}).
Therefore we can define a Lebesgue-Stieltjes integral
$\int _0^{\infty}\phi(s)dL _i(s,a,\omega)$
for each $a$ for any Borel-measurable function $\phi$. In particular
\begin{eqnarray}
\int _0^{\infty}1_{R\setminus\{a\}}(X_i(s))dL_i(s,a,\omega)=0, \  \ a.s.\ \ i=1,2.
\end{eqnarray}
Furthermore if $\phi$ is differentiable, then we have the following
integration by parts formula
\begin{eqnarray}
&&\int _0^t\phi(s)dL_i(s,a,\omega)\nonumber\\
&=&\phi(t)L_i(t,a,\omega)-\int
_0^t\phi^{\prime}(s)L_i(s,a,\omega){\rm d}s,\ \ a.s.\ \ i=1, 2.
\end{eqnarray}
Moreover, if $g(s,x_i,\omega)$ is measurable and bounded,
by the occupation times formula (e.g. see \cite{ks}, \cite{yor})),
\begin{eqnarray*}
\int _0^tg(s,X_i(s))d<\hskip-4pt M_i\hskip-4pt> _s=2\int _{-\infty}^{\infty}\int
_0^tg(s,a)dL_i(s,a,\omega){\rm d}a.\ \ a.s.\ \ i=1,2
\end{eqnarray*}
If $g(s,x_i)$ is differentiable in $s$, then using the integration by
parts formula, we have
\begin{eqnarray}
&&\int _0^tg(s,X_i(s))d<\hskip-4pt M_i\hskip-4pt> _s\nonumber\\
&=&2\int _{-\infty}^{\infty}\int
_0^tg(s,a)dL_i(s,a,\omega){\rm d}a\nonumber\\
&=&2\int _{-\infty}^{\infty}g(t,a)L_i(t,a,\omega){\rm d}a\nonumber\\
&&-2\int _{-\infty}^{\infty}\int _0^t{\partial \over \partial
s}g(s,a)L_i(s,a,\omega){\rm d}s{\rm d}a, \ \ a.s.,
\end{eqnarray}
for $i=1,2$.
On the other hand, by Tanaka formula
\begin{eqnarray*}
L_1(t,a)=(X_1(t)-a)^+-(X_1(0)-a)^+-\hat M_1(t,a)-\hat V_1(t,a),
\end{eqnarray*}
 where $ \hat Z_1(t,a)=\int_0^t1_{\{X_1(s)>a\}}dZ_1(s),\ Z_1=M_1,V_1,X_1$.
By a standard localizing argument, we may assume without loss of generality
that there is a constant $N$ for which
\begin{eqnarray*}
\sup\limits_{0\leq s\leq t} |X_1(s)|\leq N, \ <\hskip-4pt M_1\hskip-4pt>_t\leq N, \ Var_tV_1\leq N,
\end{eqnarray*}
where $ Var_tV_1$ is the total variation of $V_1$ on $[0,t]$.
From the property of local time (see Chapter $3$ in \cite {ks}), for any $\gamma\geq 1$,
 \begin{eqnarray*}
 &&E|\hat M_1(t,a)-\hat M_1(t,b)|^{2\gamma}\\
 &=&E|\int_0^t1_{\{a<X_s\leq b\}}d<\hskip-4pt M_1\hskip-4pt>_s|^\gamma\\
 &\leq&C(b-a)^\gamma, \ a<b
 \end{eqnarray*}
 where the constant $C$ depends on $\gamma$ and on the bound $N$.
 From Kolmogorov's tightness criterion (see \cite{kun}), we know that the sequence $Y_n(a):= {1\over n}\hat M_1(t,a)$, $n=1,2,\cdots$, is tight. Moreover for any $a_1, a_2,\cdots, a_k$,
 \begin{eqnarray*}
 &&P(\sup\limits_{a_i}|{1\over n}\hat M_1(t,a_i)|\leq 1)\\
 &=&P(|{1\over n}\hat M_1(t,a_1)|\leq 1, |{1\over n}\hat M_1(t,a_2)|\leq 1,\cdots, |{1\over n}\hat M_1(t,a_k)|\leq 1|)\\
 &\geq&1-\sum\limits_{i=1}^k P(|{1 \over n}\hat M_1(t,a_i)|>1)\\
 &\geq&1-{1\over n^2}\sum\limits_{i=1}^k E[\hat M_1^2(t,a_i)]\\
 &\geq&1-{k\over{n^2}}C(N-a),
 \end{eqnarray*}
 so by the weak convergence theorem of random fields (see Theorem 1.4.5 in \cite{kun}), we have
 \begin{eqnarray*}
 \lim\limits_{n\to \infty}P(\sup\limits_{a}|\hat M_1(t,a)|\leq n)=1.
 \end{eqnarray*}
 Furthermore it is easy to see that
 \begin{eqnarray*}
 {1\over n} \hat V_1(t,a)\leq {1\over n} Var_tV_1(t,a)\to 0, \ when\ n\to \infty,
 \end{eqnarray*}
so it follows that,
 \begin{eqnarray*}
  \lim\limits_{n\to \infty}P(\sup\limits_{a}|L_1(t,a)|\leq n)=1.
  \end{eqnarray*}
Therefore in our  localization argument, we can also assume $L_1(t,a)$ and
$L_2(t,a)$ are bounded uniformly in $a$.
\vskip5pt

In the following we assume some conditions on  $f:{R^+}\times R\times R \to R$:

\vskip5pt
 {\it  Condition (i)} $f(\cdot,\cdot,\cdot): {R^+}\times R\times R \to R$ is left continuous and locally bounded and jointly continuous from the right in $t$ and left in $(x_1,x_2)$ at each point $(0,x_1,x_2)$;
\vskip3pt

 {\it Condition (ii)} the left derivative ${\partial^-\over \partial t}f(t,x_1,x_2)$ exists at all points of $(0,\infty)\times R^2$, and ${\nabla}_1^-f(t,x_1,x_2)$, ${\nabla}_2^-f(t,x_1,x_2)$ exist  at all points $[0,\infty)\times R^2$ and are jointly left continuous and locally bounded;
     \vskip3pt

  {\it  Condition (iii)} ${\nabla}_i^-f(t,x_1,x_2)$ is of locally bounded variation in $x_i$, $i=1,2$;
     \vskip3pt

 {\it   Condition (iv)} ${\partial ^-\over \partial t } \nabla_i^- f(t,x_1,x_2)$ $(i=1,2)$ and $\nabla_1^-\nabla_2^- f(t,x_1,x_2)$ exist at all points of $(0,\infty)\times R^2$ and $[0,\infty)\times R^2$ respectively, and are left continuous and locally bounded;
     \vskip3pt

{ \it  Condition (v)}  $\nabla_1^-\nabla_2^- f(t,x_1,x_2)$ is of locally bounded variation in $(t,x_1)$ and $(t,x_2)$ and  $\nabla_1^-\nabla_2^- f(0,x_1,x_2)$ is of locally bounded variation in $x_1$ and $x_2$ respectively.

       \vskip 5pt
       From the assumption of $\nabla _1^-f$, we can use the one-dimensional generalized It$\hat {\rm o}$ formula (Theorem 1.1 in \cite {Zhao1})
\begin{eqnarray}
&&\nabla _1^-f(t,a,X_2(t))-\nabla _1^-f(0,a,X_2(0))\nonumber\\
&=&\int _0^t{\partial ^-\over \partial s} \nabla _1^-f(s,a,X_2(s)){\rm d}s+ \int _0^t\nabla _1^-\nabla _2 ^-f(s,a,X_2(s))dX_2(s)\nonumber\\
&&+
\int _{-\infty}^{\infty }{L}_2(t,x_2){\rm d}_{x_2}\nabla _1^-\nabla _2 ^-f(t,a,x_2)\nonumber\\
&&-\int _{-\infty}^{+\infty}\int _0^{t}{L}_2(s,x_2){\bf \rm d}_{s,x_2}\nabla _1^-\nabla _2^-f(s,a,x_2).\ \ a.s.
\end{eqnarray}
Therefore $\nabla _1^-f(t,a,X_2(t))$ is a continuous
semi-martingale, and can be decomposed as $\nabla
_1^-f(t,a,X_2(t))=\nabla _1^-f(0,a,X_2(0))+ h(t,a) + v(t,a)$,
where $h$ is a continuous local martingale and $v$ is a continuous
process of locally bounded variation (in $t$). In fact
$h(t,a)=\int _0^t\nabla _1^-\nabla _2 ^-f(s,a,X_2(s))dM_2(s)$.
  Define
\begin{eqnarray*}
F_s(a,b)&:=&<h(a),h(b)>_s\ =\
<\nabla_1^-f(a),\nabla_1^-f(b)>_s\nonumber\\&=&\int _0^s\nabla
_1^-\nabla _2 ^-f(r,a,X_2(r))\nabla _1^-\nabla _2
^-f(r,b,X_2(r)){\bf \rm d}<\hskip-4pt M_2\hskip-4pt>_r.\\
F_{s_{k}}^{s_{k+1}}(a,b)&:=&<h(a),h(b)>_{s_{k}}^{s_{k+1}}\ =\
<\nabla_1^-f(a),\nabla_1^-f(b)>_{s_{k}}^{s_{k+1}}\nonumber\\&=&\int
_{s_{k}}^{s_{k+1}}\nabla _1^-\nabla _2 ^-f(r,a,X_2(r))\nabla
_1^-\nabla _2 ^-f(r,b,X_2(r)){\bf \rm d}<\hskip-4pt
M_2\hskip-4pt>_r.
\end{eqnarray*}
  We need to prove $h(s,a) \in {\cal V}_2$.
To see this, as $\nabla _1^-\nabla _2 ^-f(t,x_1,x_2)$ is of
locally bounded variation in $x_1$, so for any compact set $G$,
$\nabla _1^-\nabla _2 ^-f(t,x_1,x_2)$ is of bounded variation in
$x_1$ for $x_1 \in G$. Also on this set, let $\cal P$ be the
partition on $R^2\times [0,t]$, ${\cal P}_i$ be a partition on $R$
$(i=1,2)$, ${\cal P}_3$ be a partition on $[0,t]$  such that
${\cal P} = {{\cal P}_1}\times {\cal P}_2\times {\cal P}_3$. Then
we have:
\begin{eqnarray*}
&&{\rm Var} _{s,a,b} (F_{s}(a,b))\\&=&\sup_{\cal
P}\sum_k\sum_i\sum_j\Big|F_{s_{k}}^{s_{k+1}}(a_{i+1},b_{j+1}) -
F_{s_{k}}^{s_{k+1}}(a_{i+1},b_{j}) -
F_{s_{k}}^{s_{k+1}}(a_{i},b_{j+1})\\
&&\hskip 2.5cm+F_{s_{k}}^{s_{k+1}}(a_{i},b_{j})\Big|\\
&=&\sup _{\cal P}\sum_k\sum_i\sum_j\Big|\int_{s_{k}}^{s_{k+1}}\nabla _1^-\nabla _2^-f(r,a_{i+1},X_2(r))\nabla _1^-\nabla _2^-f(r,b_{j+1},X_2(r)){\bf \rm d}<\hskip-4pt M_2\hskip-4pt>_r \\
&&- \int_{s_{k}}^{s_{k+1}}\nabla _1^-\nabla _2^-f(r,a_{i+1},X_2(r))\nabla _1^-\nabla _2^-f(r,b_{j},X_2(r)){\bf \rm d}<\hskip-4pt M_2\hskip-4pt>_r \\
&&- \int_{s_{k}}^{s_{k+1}}\nabla _1^-\nabla _2^-f(r,a_{i},X_2(r))\nabla _1^-\nabla _2^-f(r,b_{j+1},X_2(r)){\bf \rm d}<\hskip-4pt M_2\hskip-4pt>_r\\
&&+ \int_{s_{k}}^{s_{k+1}}\nabla _1^-\nabla _2^-f(r,a_{i},X_2(r))\nabla _1^-\nabla _2^-f(r,b_{j},X_2(r)){\bf \rm d}<\hskip-4pt M_2\hskip-4pt>_r  \bigg| \\
&=&\sup _{\cal
P}\sum_k\sum_i\sum_j\bigg|\int_{s_{k}}^{s_{k+1}}\bigg(\nabla
_1^-\nabla _2^-f(r,a_{i+1},X_2(r))
- \nabla _1^-\nabla _2^-f(r,a_{i},X_2(r))\bigg)\\
&&\bigg( \nabla _1^-\nabla _2^-f(r,b_{j+1},X_2(r)) - \nabla _1^-\nabla _2^-f(r,b_{j},X_2(r))\bigg){\bf \rm d}<\hskip-4pt M_2\hskip-4pt>_r\bigg |\\
&\leq& \int_0^{s} \sup _{{\cal P}_1}\sum_i\Big|\nabla _1^-\nabla _2^-f(r,a_{i+1},X_2(r)) - \nabla _1^-\nabla _2^-f(r,a_{i},X_2(r))\Big|\\
&&\sup_{{\cal P}_2}\sum_j\Big| \nabla _1^-\nabla _2^-f(r,b_{j+1},X_2(r))- \nabla _1^-\nabla _2^-f(r,b_{j},X_2(r))\Big|{\bf \rm d}<\hskip-4pt M_2\hskip-4pt>_r \\
&=&\int_0^{s}\bigg({\rm Var}_{a}(\nabla _1^-\nabla _2^-f(r,a,X_2(r)))\bigg)^2{\bf \rm d}<\hskip-4pt M_2\hskip-4pt>_r <{\infty}.
\end{eqnarray*}
Therefore under the localizing assumption, $\int_{-\infty}^\infty\int_0^t L_1(s,a)d_{s,a} \nabla_1^-f(s,a,X_2(s))$ and $\int_{-\infty}^\infty\int_0^t L_2(s,a)d_{s,a} \nabla_2^-f(s,X_1(s),a)$ can be defined by Definition \ref{definition1}. A localizing argument implies they are semi-martingales.

    We will prove the following generalized It$\hat{\rm o}$'s formula in two-dimensional space.

\begin{thm} \label{tom100}
Under  conditions (i)-(v), for any continuous two-dimensional semi-martingale $X(t)=(X_1(t), X_2(t))$, we have
\begin{eqnarray}
&&
f(t,X_1(t),X_2(t))\nonumber\\
&=&f(0,X_1(0),X_2(0))+\int _0^t{\partial ^-\over \partial s} f(s,X_1(s),X_2(s)){\rm d}s\nonumber\\
&&+\sum_{i=1}^2\int _0^t\nabla _i ^-f(s,X_1(s),X_2(s))dX_i(s)\nonumber\\
&&+\int _{-\infty}^{\infty }L _1(t,a){\rm d}_a\nabla _1 ^-f(t,a,X_2(t))
-\int _{-\infty}^{+\infty}\int _0^{t}L _1(s,a)
{\bf \rm d}_{s,a}\nabla _1^-f(s,a,X_2(s))\nonumber\\
&&+\int _{-\infty}^{\infty }L _2(t,a){\rm d}_a\nabla _2 ^-f(t,X_1(t),a)
-\int _{-\infty}^{+\infty}\int _0^{t}L _2(s,a)
{\bf \rm d}_{s,a}\nabla _2^-f(s,X_1(s),a)\nonumber\\
&&+\int_0^t{\nabla _1^-}{\nabla _2^-}f(s,X_1(s),X_2(s))d<\hskip-4pt M_1,M_2\hskip-4pt>_s.\ \ a.s.
\end{eqnarray}
\end{thm}
{\em Proof}: By a standard localization argument, we can assume $X_1(t)$, $X_2(t)$ and their quadratic variations $<\hskip-4pt X_1\hskip-4pt>_t$,$<\hskip-4pt X_2\hskip-4pt>_t$ and $<\hskip-4pt X_1,X_2\hskip-4pt>_t$ and the local times $L_1$, $L_2$ are bounded processes so that $f$, ${\partial ^-\over \partial t} f $, $\nabla_1^-f$, $\nabla_2^-f$, $\nabla_1^-\nabla_2^-f$, $Var_{x_1}\nabla_1^-f$, $Var_{x_2}\nabla_2^-f$, $Var_{x_1}\nabla_1^-\nabla_2^-f$, $Var_{x_2}\nabla_1^-\nabla_2^-f$ are bounded. Note the left derivatives of $f$ agree with the generalized derivatives, so condition (ii) and (iv) imply that $f$ is absolutely continuous in each variable and $\nabla_1^-f$ is absolutely continuous with respect to $t$ and $x_2$ respectively and $\nabla_2^-f$ is absolutely continuous with respect to $t$ and $x_1$ respectively.

We divide the proof into several steps:

 {\bf (A)} Define
\begin{eqnarray}\label{smooth}
\rho(x)=\cases {c{\rm e}^{{1\over (x-1)^2-1}}, {\rm \  if } \ x\in (0,2),\cr
0, \ \ \ \ \ \ \ \ \ \ \ \ \  {\rm otherwise.}}
\end{eqnarray}
Here $c$ is chosen such that $\int _0^2\rho(x)dx=1$.
Take $\rho_n(x)=n\rho(nx)$ as mollifiers. Define
\begin{eqnarray*}
&&f_n(s,x_1,x_2)\\
&=&\int _{-\infty}^{+\infty}\int _{-\infty}^{+\infty}\int
_{-\infty}^{+\infty}\rho_n(s-\tau)\rho_n(x_1-y)\rho_n(x_2-z)f(\tau,y,z)d\tau dydz, \
n\geq 1,
\end{eqnarray*}
where we set $f(\tau,y,z)=f(-\tau,y,z)$ if $\tau<0$.
Then $f_n(s,x_1,x_2)$ are smooth and
\begin{eqnarray}\label{truman41}
&&f_n(s,x_1,x_2)\nonumber\\
&=&\int _0^2\int _0^2\int _0^2\rho(t)\rho(y)\rho(z)f(s-{t\over
n},x_1-{y\over n},x_2-{z\over n})dtdydz, \ n\geq 1.
\end{eqnarray}
Because of the absolute continuity mentioned above, we can differentiate under the integral (\ref {truman41}) to see ${\partial \over \partial t}f_n$, $\nabla_1f_n$, $\nabla_2f_n$, $\nabla_1\nabla_2f_n$, $Var_{x_1}\nabla_1f_n$, $Var_{x_2}\nabla_2f_n$, $Var_{x_1}\nabla_1\nabla_2f_n$ and $Var_{x_2}\nabla_1\nabla_2f_n$ are bounded. Furthermore using Lebesgue's dominated convergence theorem, one can prove that as
$n\to \infty$,
\begin{eqnarray}
f_n(s,x_1,x_2)&\to & f(s,x_1,x_2),\ s\geq 0{\label{frc1}}\\
{\partial \over \partial s}f_n(s,x_1,x_2)&\to & {\partial ^- \over \partial s} f(s,x_1,x_2),\ s> 0{\label{frc2}}\\
\nabla_1f_n(s,x_1,x_2)&\to&\nabla _1^-f(s,x_1,x_2),\ s\geq 0{\label{frc3}}\\
\nabla_2f_n(s,x_1,x_2)&\to & \nabla _2^-f(s,x_1,x_2),\ s\geq 0{\label{frc4}}\\
\nabla_1\nabla _2f_n(s,x_1,x_2)&\to&\nabla_1^-\nabla _2^-f(s,x_1,x_2),\ s\geq 0{\label{frc5}}
\end{eqnarray}
and each $(x_1,x_2) \in R^2$.

\vskip5pt
\noindent
{\bf (B)}
It turns out for any $g(t,x_1)$ being continuous in $t$ and $C^1$ in $x_1$ and having a compact support, using the integration by parts
formula and Lebesgue's
dominated convergence theorem, we see that
\begin{eqnarray}
&&\lim_{n\to +\infty}\int _{-\infty}^{+\infty} g(t,x_1)\Delta_1 f_n(t,x_1,X_2(t))dx_1\nonumber\\
&=&
-\lim_{n\to +\infty}\int _{-\infty}^{\infty}
\nabla  g(t,x_1)\nabla _1f_n(t,x_1,X_2(t))dx_1
\nonumber\\
&=&
-\int _{-\infty}^{\infty}\nabla  g(t,x_1)\nabla _1^-f(t,x_1,X_2(t))dx_1.
\end{eqnarray}
Note $\nabla _1^-f(t,x_1,x_2)$ is of locally bounded variation in $x_1$ and $g(t,x_1)$
has a compact support in $ x_1$, so
\begin{eqnarray}
&&-\int _{-\infty}^{+\infty} \nabla  g(t,x_1)\nabla _1^-
f(t,x_1,X_2(t))dx_1\nonumber\\
&=&\int _{-\infty}^{
+\infty} g(t,x_1){\rm d}_{x_1}\nabla _1^-f(t,x_1,X_2(t)).
\end{eqnarray}
Thus
\begin{eqnarray}\label{fcr2}
&&\lim_{n\to +\infty}\int _{-\infty}^{+\infty} g(t,x_1)\Delta _1f_n(t,x_1,X_2(t))dx_1\nonumber\\
&=&\int _{-\infty}^{\infty } g(t,x_1){\rm d}_{x_1}\nabla _1^-f(t,x_1,X_2(t)).
\end{eqnarray}
\vskip 5pt
\noindent
{\bf (C)}
If $g(s,x_1)$ is $C^2$ in $x_1$, $\Delta g(s,x_1)$ is bounded uniformly in $s$, ${\partial \over \partial s}\nabla g(s,x_1)$ is
continuous in $s$ and has a compact support in $x_1$, and \\
$E\left[\int _0^t\int_{R^2}|g(s,x)g(s,y)||{\bf \rm
d}_{x,y,s}<h(x),h(y)>_s|\right] < \infty$, where $h\in {\cal
V}_2$, then applying It$\hat{\rm o}$'s formula, Lebesgue's
dominated convergence theorem and the integration by parts
formula,
\begin{eqnarray*}
&&\lim_{n\to +\infty}\Big(\int _{0}^{t}\int _{-\infty}^{+\infty}g(s,x_1){\partial \over \partial s}\Delta _1f_n(s,x_1,X_2(s))dx_1{\rm d}s \\
&&+ \int _{0}^{t}\int _{-\infty}^{+\infty}g(s,x_1)\nabla _2\Delta _1f_n(s,x_1,X_2(s))dx_1dX_2(s) \\
&&+{1\over2} \int _{0}^{t}\int _{- \infty}^{+\infty}g(s,x_1)\Delta_2\Delta _1f_n(s,x_1,X_2(s))dx_1d<\hskip-4pt M_2\hskip-4pt>_s\Big) \\
&=&
-\lim_{n\to +\infty}\Big(\int _{0}^{t}\int _{-\infty}^{+\infty}
\nabla g(s,x_1){\partial \over \partial
s}\nabla _1f_n(s,x_1,X_2(s))dx_1{\rm d}s \\
&&+ \int _{0}^{t}\int _{-\infty}^{+\infty}\nabla g(s,x_1)\nabla _1\nabla_2f_n(s,x_1,X_2(s))dx_1dX_2(s) \\
&&+{1\over2}\int _{0}^{t}\int _{- \infty}^{+\infty}\nabla g(s,x_1)\Delta _2\nabla_1f_n(s,x_1,X_2(s))dx_1d<\hskip-4pt M_2\hskip-4pt>_s \Big)\\
&=&-\lim_{n\to +\infty}\int _{-\infty}^{\infty}\int _0^t \nabla g(s,x_1){\rm d}_s\nabla _1f_n(s,x_1,X_2(s))dx_1\\
&=&-\lim_{n\to +\infty}\Big(\int _{-\infty}^{\infty} \nabla g(s,x_1)\nabla _1f_n(s,x_1,X_2(s))\Big |_0^t dx_1\\
&&-\int _{0}^{t}\int _{-\infty}^{+\infty}
{\partial \over \partial s}\nabla g(s,x_1)\nabla _1f_n(s,x_1,X_2(s))dx_1{\rm d}s\Big)\\
&=&-\int _{-\infty}^{\infty} \nabla g(s,x_1)\nabla _1^-f(s,x_1,X_2(s))\Big |_0^t dx_1\\
&&+\int _{0}^{t}\int _{-\infty}^{+\infty}
{\partial \over \partial s}\nabla g(s,x_1)\nabla _1^-f(s,x_1,X_2(s))dx_1{\rm d}s\\
&=&-\int _{-\infty}^{+\infty}\int _0^t \nabla g(s,x_1) {\rm d}_{s} \nabla _1^-f(s,x_1,X_2(s))dx_1.
\end{eqnarray*}
It turns out by applying Proposition \ref{proposition1} that
\begin{eqnarray}\label {fcr3}
&&\lim_{n\to +\infty}\Big(\int _{0}^{t}\int _{-\infty}^{+\infty}g(s,x_1){\partial \over \partial s}\Delta _1f_n(s,x_1,X_2(s))dx_1{\rm d}s \nonumber\\
&&+ \int _{0}^{t}\int _{-\infty}^{+\infty}g(s,x_1)\nabla _2\Delta _1f_n(s,x_1,X_2(s))dx_1dX_2(s) \nonumber\\
&&+{1\over2} \int _{0}^{t}\int _{- \infty}^{+\infty}g(s,x_1)\Delta_2\Delta _1f_n(s,x_1,X_2(s))dx_1d<\hskip-4pt M_2\hskip-4pt>_s\Big) \nonumber\\
&=&\int _0^t \int _{-\infty}^{+\infty}g(s,x_1){\rm \bf d}_{s,x_1} \nabla _1^-f(s,x_1,X_2(s)).
\end{eqnarray}
\vskip5pt
\noindent
{\bf (D)}
But any c$\grave{a}$dl$\grave{a}$g function with
a compact support can be approximated by
smooth
functions with a compact support uniformly
by the following
standard smoothing procedure
\begin{eqnarray*}
g_m(t,x_1)=\int
_{-\infty}^{\infty}\rho_m(y-x_1)g(t,y)dy=\int _0^2\rho (z)g(t,x_1+{z\over
m})dz.
\end{eqnarray*}
Then we can prove that (\ref{fcr2}) also holds for any c$\grave{a}$dl$\grave{a}$g function $g(t,x_1)$ with a
compact support in $x_1$. Moreover, if $g \in {\cal V}_3$, (\ref{fcr3}) also holds.

To see (\ref{fcr2}), note that there is a compact set $G\subset R^1$ such that
\begin{eqnarray*}
&\max \limits _{x_1\in G}|g_m(t,x_1)-g(t,x_1)|\to 0 & {\rm as } \ \ m\to +\infty,\\
&g_m(t,x_1)=g(t,x_1) =0 & {\rm for } \ \ x_1\notin G.
\end{eqnarray*}
Note
\begin{eqnarray}\label{elworthy1}
&&\int _{-\infty}^{+\infty}g(t,x_1)\Delta_1 f_n(t,x_1,X_2(t))dx_1\nonumber\\
&=&
\int _{-\infty}^{+\infty}g_m(t,x_1)\Delta_1 f_n(t,x_1,X_2(t))dx_1\nonumber\\
&&+\int _{-\infty}^{+\infty}(g(t,x_1)-g_m(t,x_1))\Delta_1 f_n(t,x_1,X_2(t))dx_1.
\end{eqnarray}
It is easy to see from ({\ref{fcr2}) and Lebesgue's
dominated convergence theorem, that
\begin{eqnarray}\label{elworthy3}
&&
\lim\limits_{m\to \infty}\lim\limits _{n\to \infty}
\int _{-\infty}^{+\infty}g_m(t,x_1)\Delta_1 f_n(t,x_1,X_2(t))dx_1\nonumber\\
&=&\lim\limits _{m\to \infty}\int _{-\infty}^{\infty } g_m(t,x_1){\bf \rm
d}_{x_1}\nabla_1^-f(t,x_1,X_2(t))\nonumber\\
&=&\int _{-\infty}^{\infty } g(t,x_1){\bf \rm d}_{x_1}\nabla_1^-f(t,x_1,X_2(t)).
\end{eqnarray}
Moreover,
\begin{eqnarray}\label{elworthy6}
&&|\int _{-\infty}^{+\infty}\Big(g(t,x_1)-g_m(t,x_1)\Big)\Delta _1f_n(t,x_1,X_2(t))dx_1|\nonumber\\
&=& |\int _{-\infty}^{+\infty}\Big(g(t,x_1)-g_m(t,x_1)\Big)d_{x_1}\nabla_1f_n(t,x_1,X_2(t))|\nonumber\\
&\leq&
\Big(\max\limits _{x_1\in {G}} |g(t,x_1)-g_m(t,x_1)|\Big)Var_{x_1\in G}\nabla _1f_n(t,x_1,X_2(t)).
\end{eqnarray}
But,
\begin{eqnarray*}
\lim\limits_{m\to \infty}\limsup\limits_{n\to \infty}
\Big(\max\limits _{x_1\in {G}} |g(t,x_1)-g_m(t,x_1)|\Big)Var_{x_1\in G}\nabla _1f_n(t,x_1,X_2(t))=0.
\end{eqnarray*}
So inequality (\ref{elworthy6}) leads to
\begin{eqnarray}\label{elworthy7}
&&\lim\limits_{m\to \infty}\limsup\limits_{n\to \infty}
|\int _{-\infty}^{+\infty}\Big(g(t,x_1)-g_m(t,x_1)\Big)\Delta _1f_n(t,x_1,X_2(t))dx_1|\nonumber\\
&=&0.
\end{eqnarray}
Now we use (\ref{elworthy1}), (\ref{elworthy3}) and (\ref{elworthy7})
\begin{eqnarray*}
&&
\limsup\limits_{n\to \infty}\int _{-\infty}^{+\infty}g(t,x_1)\Delta_1
f_n(t,x_1,X_2(t))dx_1\nonumber\\
&=&
\lim\limits_{m\to \infty}\limsup\limits_{n\to \infty}\int
_{-\infty}^{+\infty}g_m(t,x_1)\Delta _1f_n(t,x_1,X_2(t))dx_1\nonumber\\
&&
+\lim\limits_{m\to \infty}\limsup\limits_{n\to \infty}\int
_{-\infty}^{+\infty}\Big(g(t,x_1)-g_m(t,x_1)\Big)\Delta _1f_n(t,x_1,X_2(t))dx_1\nonumber
\\
&=& \int _{-\infty}^{\infty } g(t,x_1){\rm }{\rm d}_{x_1}\nabla _1^-f(t,x_1,X_2(t)).
\end{eqnarray*}
Similarly we also have
\begin{eqnarray}\label{elworthy7*}
&&\liminf\limits_{n\to \infty}\int _{-\infty}^{+\infty}g(t,x_1)\Delta _1f_n(t,x_1,X_2(t))dx_1\nonumber\\
&=&
\int _{-\infty}^{\infty } g(t,x_1){\rm d}_{x_1}\nabla _1^-f(t,x_1,X_2(t)).
\end{eqnarray}
So (\ref{fcr2}) holds for a c$\grave{a}$dl$\grave{a}$g function $g$ with a compact support in $x_1$.

Now we prove that (\ref{fcr3}) also holds for a c$\grave{a}$dl$\grave{a}$g function $g\in {\cal V}_3$.
 Obviously,
 \begin{eqnarray*}
&&\int _{0}^{t}\int _{-\infty}^{+\infty}g(s,x_1){\partial \over \partial s}\Delta _1f_n(s,x_1,X_2(s))dx_1{\rm d}s \\
&&+ \int _{0}^{t}\int _{-\infty}^{+\infty}g(s,x_1)\nabla _2\Delta _1f_n(s,x_1,X_2(s))dx_1dX_2(s) \\
&&+{1\over2} \int _{0}^{t}\int _{- \infty}^{+\infty}g(s,x_1)\Delta_2\Delta _1f_n(s,x_1,X_2(s))dx_1d<\hskip-4pt M_2\hskip-4pt>_s \\
&=&\int _0^t \int _{-\infty}^{+\infty}g(s,x_1){\rm \bf d}_{s,x_1} \nabla _1f_n(s,x_1,X_2(s)).
\end{eqnarray*}
Define
\begin{eqnarray*}
g_m(s,x_1)=\int_{-\infty}^{\infty}\int_{-\infty}^{\infty}\rho_m(y-x_1)\rho_m(\tau-s)g(\tau,y)d\tau dy.
\end{eqnarray*}
Then there is a compact $G\subset R^1$ such that
\begin{eqnarray*}
&\max \limits _{0\leq s \leq t, x_1\in G}|g_m(s,x_1)-g(s,x_1)|\to 0 & {\rm as } \ \ m\to +\infty,\\
&g_m(s,x_1)=g(s,x_1) =0 & {\rm for } \ \ x_1\notin G.
\end{eqnarray*}
Then it is trivial to see
\begin{eqnarray*}
&&\int _0^t \int _{-\infty}^{+\infty}g(s,x_1){\rm \bf d}_{s,x_1} \nabla _1f_n(s,x_1,X_2(s))\\
&=&\int _0^t \int _{-\infty}^{+\infty}g_m(s,x_1){\rm \bf d}_{s,x_1} \nabla _1f_n(s,x_1,X_2(s))\\
&&+\int _0^t \int _{-\infty}^{+\infty}(g(s,x_1)-g_m(s,x_1)){\rm \bf d}_{s,x_1} \nabla _1f_n(s,x_1,X_2(s)).
\end{eqnarray*}
But from (\ref{fcr3}), we can see that
\begin{eqnarray}
&&\lim\limits_{m\to\infty}\lim\limits_{n\to\infty}\int _0^t \int _{-\infty}^{+\infty}g_m(s,x_1){\rm \bf d}_{s,x_1} \nabla _1f_n(s,x_1,X_2(s))\nonumber\\
&=&\lim\limits_{m\to\infty}\int _0^t \int _{-\infty}^{+\infty}g_m(s,x_1){\rm \bf d}_{s,x_1} \nabla _1^-f(s,x_1,X_2(s))\nonumber\\
&=&\int _0^t \int _{-\infty}^{+\infty}g(s,x_1){\rm \bf d}_{s,x_1} \nabla _1f(s,x_1,X_2(s)). \ \ \ \ \     (limit\  in\  {\cal M}_2)
\end{eqnarray}
The last limit holds because of the following:
\begin{eqnarray*}
&&E\Big[\int _{-\infty}^{+\infty}(g_m(s,x_1)-g(s,x_1)){\rm \bf d}_{s,x_1} \nabla _1^-f(s,x_1,X_2(s))\Big]^2\\
&=&E\Big[\int_0^t\int_{-\infty}^{+\infty}(g_m-g)(s,a)(g_m-g)(s,b) {\rm\bf d}_{a,b} {\rm d}_{s} <\nabla _1^-f(a),\nabla _1^-f(b)>_s\Big]\\
&=&E\Big[\int_0^t\int_{-\infty}^{+\infty}(g_m-g)(s,a)(g_m-g)(s,b)\\
&&\ \ \ \ \ {\rm \bf d}_{a,b} \nabla_1^-\nabla _2^-f(s,a,X_2(s))\nabla_1^-\nabla _2^-f(s,b,X_2(s))\Big]
{\rm d}<\hskip-4pt M_2\hskip-4pt>_s\\
&=&\int_0^tE\Big[\int_{-\infty}^{+\infty}(g_m-g)(s,a) {\rm d}_{a} \nabla_1^-\nabla _2^-f(s,a,X_2(s))\Big]^2 {\rm d}<\hskip-4pt M_2\hskip-4pt>_s \\
&\to& 0,\ as  \ m\to {\infty}.\\
\end{eqnarray*}
On the other hand,
\begin{eqnarray}\label{elworthy7**}
&&\lim\limits_{m\to\infty}\lim\limits_{n\to\infty}\int _0^t \int _{-\infty}^{+\infty}(g(s,x_1)-g_m(s,x_1)){\rm \bf d}_{s,x_1} \nabla _1f_n(s,x_1,X_2(s))  \nonumber\\
&=&0. \ \ \ \ \ \ \ \  \ \     (limit\  in\ {\cal M}_2)
\end{eqnarray}
In fact,
\begin{eqnarray*}
&&E\Big[\int _{-\infty}^{+\infty}(g(s,x_1)-g_m(s,x_1)){\rm \bf d}_{s,x_1} \nabla _1^-f_n(s,x_1,X_2(s))\Big]^2\\
&=&\int_0^tE\Big[\int_{-\infty}^{+\infty}(g-g_m)(s,a) {\rm d}_{a} \nabla_1\nabla _2f_n(s,a,X_2(s))\Big]^2 d<\hskip-4pt M_2\hskip-4pt>_s.
\end{eqnarray*}
Note that $\nabla_1\nabla _2f_n(s,a,X_2(s))$ is of bounded variation in $a$, we can use an argument similar to the one in the proof of (\ref{elworthy7}) and (\ref{elworthy7*}) to prove (\ref{elworthy7**}).
\vskip5pt
\noindent
{\bf (E)}
Now we use the multi-dimensional It${\hat {\rm o}}$'s formula to the function \\
$f_n(s,X_1(s),X_2(s))$, then a.s.
\begin{eqnarray}\label{zhao11}
&&f_n(t,X_1(t),X_2(t))-f_n(0,X_1(0),X_2(0))\nonumber\\&=&\int _0^t{\partial \over \partial s}
f_n(s,X_1(s),X_2(s)){\rm d}s
+\sum_{i=1}^{2}\int _0^t\nabla_i f_n(s,X_1(s),X_2(s))dX_i(s)\nonumber\\
&&+{1\over 2}\int _0^t\Delta _1f_n(s,X_1(s),X_2(s))d<\hskip-4pt M_1\hskip-4pt>_s\nonumber\\
&&+
{1\over 2}\int _0^t\Delta _2f_n(s,X_1(s),X_2(s))d<\hskip-4pt M_2\hskip-4pt>_s\nonumber\\
&&+\int _0^t\nabla _1\nabla _2f_n(s,X_1(s),X_2(s))d<\hskip-4pt M_1,M_2\hskip-4pt>_s.
\end{eqnarray}
As $n\to \infty$, it is easy to see from Lebesgue's dominated
convergence theorem and (\ref{frc1}), (\ref{frc2}), (\ref{frc3}), (\ref{frc4}), (\ref{frc5}) that, $(i=1,2)$
\begin{eqnarray*}
f_n(t,X_1(t),X_2(t))-f_n(0,z_1,z_2)\to f(t,X_1(t),X_2(t))-f(0,z_1,z_2),\ \ a.s.
\end{eqnarray*}
\begin{eqnarray*}
\int _0^t{\partial \over \partial s} f_n(s,X_1(s),X_2(s)){\rm d}s\to \int
_0^t{\partial ^- \over \partial s} f(s,X_1(s),X_2(s)){\rm d}s,\ \ a.s.
\end{eqnarray*}
\begin{eqnarray*}
\int _0^t \nabla _if_n(s,X_1(s),X_2(s))dV_i(s)\to \int _0^t\nabla _i^-f(s,X_1(s),X_2(s))dV_i(s),\ \ a.s
\end{eqnarray*}
\begin{eqnarray*}
&&
\int _0^t \nabla _1\nabla _2f_n(s,X_1(s),X_2(s))d<\hskip-4pt M_1,M_2\hskip-4pt>_s\\
&\to& \int _0^t\nabla _1^-\nabla _2^- f(s,X_1(s),X_2(s))d<\hskip-4pt M_1,M_2\hskip-4pt>_s.\ \ a.s.
\end{eqnarray*}
and
\begin{eqnarray*}
&&E\int _0^t(\nabla _if_n(s,X_1(s),X_2(s)))^2d<\hskip-4pt M_i\hskip-4pt>_s\\
&\to& E \int _0^t(\nabla _i
^-f(s,X_1(s),X_2(s))^2d<\hskip-4pt M_i\hskip-4pt>_s.
\end{eqnarray*}
Therefore in ${\cal M}_2$,
\begin{eqnarray*}
\int _0^t \nabla _if_n(s,X_1(s),X_2(s))dM_i(s)\to \int _0^t\nabla _i^-f(s,X_1(s),X_2(s))dM_i(s), (i=1,2).
\end{eqnarray*}
To see the convergence of ${1\over 2}\int _0^{t}\Delta _1f_n(s,X_1(s),X_2(s))d<\hskip-4pt M_1\hskip-4pt>_s$, we recall the well-known
result that the local time
$L _1(s,a)$ is jointly continuous in $s$
and c$\grave{a}$dl$\grave{a}$g with respect to $a$ and
has a compact support in space $a$ for each $s$ (\cite{yor}, \cite{ks}). As $L _1(s,a)$ is
an increasing
function of $s$ for each $a$, so if $G\subset R^1$ is
the support of $L _1(s,a)$, then
$L _1(s,a)=0$ for all $a\notin G$
and $s\le t$.  Now we use the occupation times formula, the
integration by parts formula
and (\ref{fcr2}), (\ref{fcr3}) for the case when $g$ is c$\grave{a}$dl$\grave{a}$g
with compact support,
\begin{eqnarray*}
&&
{1\over 2}\int _0^{t}\Delta _1f_n(s,X_1(s),X_2(s))d<\hskip-4pt M_1\hskip-4pt>_s\\
 &=&\int
_{-\infty}^{+\infty}\int _0^{t}\Delta _1f_n(s,a,X_2(s)){\rm d}_sL
_1(s,a){\rm d}a\nonumber\\
&=&
\int _{-\infty}^{+\infty}\Delta _1f_n({t},a,X_2({t})){L_1(t,a)}{\rm d}a\\
&&-\int
_{-\infty}^{+\infty}\bigg[\int _0^{t}{{\rm d}\over {\rm d}s}
\Delta _1f_n(s,a,X_2(s))L _1(s,a){\rm d}s\\
&&+\int _0^{t}\nabla _2\Delta _1f_n(s,X_1(s),a)L  _1(s,a)dX_2(s)\\
&&+ {1\over 2}\int _0^{t}\Delta _2\Delta _1f_n(s,X_1(s),a)L _1(s,a)d<\hskip-4pt M_2\hskip-4pt>_s\bigg]{\rm d}a\nonumber\\
&\to &
\int _{-\infty}^{\infty }L _1({t},a){\rm d}_a\nabla _1^-f({t},a,X_2({t}))\\
&&-\int
_{-\infty}^{+\infty}\int _0^{t}L _1(s,a)
{\bf \rm d}_{s,a}\nabla _1^-f(s,a,X_2(s))
\end{eqnarray*}
 as $n\to \infty$.
 About the term ${1\over 2}\int _0^{t}\Delta _2f_n(s,X_1(s),X_2(s))d<\hskip-4pt M_2\hskip-4pt>_s$, we can use the same method to get a similar result. So we proved the desired formula.
 $\hfill\diamond$\\
 \vskip 5pt
 The above smoothing procedure can be used to prove that if $f:R^+\times R^2\to R$ is left continuous and locally bounded, $C^1$ in $x_1$ and $x_2$, and the left derivatives ${\partial^-\over \partial t}f(t,x_1,x_2)$, ${{\partial^{2-}}\over {\partial{x_i}\partial{x_j}}} f(t,x_1,x_2)$, $(i,j=1,2)$ exist at all points of $(0,\infty)\times R^2$ and $[0,\infty)\times R^2$ respectively and are locally bounded and left continuous, then
 \begin{eqnarray}\label{cfr1}
&& f(t,X(t))-f(0,X(0))\nonumber\\
 &=&\int_0^t{\partial^-\over \partial s}f(s,X_1(s),X_2(s))ds
 +\sum_{i=1}^2\int_0^t\nabla_i f(s,X_1(s),X_2(s))dX_i(s)\nonumber\\
 &&+{1\over 2}\sum_{i,j=1}^2 \int_0^t {{\partial^{2-}}\over {\partial{x_i}\partial{x_j}}} f(s,X_1(s),X_2(s))d<\hskip-4pt X_i,X_j\hskip-4pt>_s.
 \end{eqnarray}
 This can be seen from the convergence in the proof of Theorem \ref{tom100} and the fact that ${{\partial^{2}}\over {\partial{x_i}\partial{x_j}}} f_n(s,x_1,x_2) \to {{\partial^{2-}}\over {\partial{x_i}\partial{x_j}}} f(s,x_1,x_2)$ under the stronger condition on ${{\partial^{2-}}\over {\partial{x_i}\partial{x_j}}}f$.

 The next theorem is an easy consequence of the methods of the proofs of Theorem \ref{tom100} and (\ref{cfr1}).

 \begin{thm}\label{cfr11}
Let $f:R^+ \times R^2\to R$ satisfy conditions (i),(ii) and $f(t,x_1,x_2)=f_h(t,x_1,x_2)+f_v(t,x_1,x_2)$.  Assume $f_h$ is $C^1$ in $x_1,x_2$ and the left derivatives ${{\partial^{2-}}\over {\partial{x_i}\partial{x_j}}}f_h(s,x_1,x_2)(i,j=1,2)$ exist and are left continuous and locally bounded; $f_v$ satisfies conditions (iii)-(v). Then
\begin{eqnarray}\label{f6}
&&f(t,X_1(t),X_2(t))-f(0,X_1(0),X_2(0))\nonumber\\
&=&\int _0^t{\partial ^-\over \partial s} f(s,X_1(s),X_2(s)){\rm d}s
+\sum_{i=1}^2\int _0^t\nabla _i ^-f(s,X_1(s),X_2(s))dX_i(s)\nonumber\\
&&+{1\over 2}\sum_{i=1}^2\int_0^t \Delta_i ^-f_h(s,X_1(s),X_2(s))d<\hskip-4pt X_i
\hskip-4pt>_s\nonumber\\
&&+\int _{-\infty}^{\infty }L _1(t,a){\rm d}_a\nabla _1 ^-f_v(t,a,X_2(t))
-\int _{-\infty}^{+\infty}\int _0^{t}L _1(s,a)
{\bf \rm d}_{s,a}\nabla _1^-f_v(s,a,X_2(s))\nonumber\\
&&+\int _{-\infty}^{\infty }L _2(t,a){\rm d}_a\nabla _2 ^-f_v(t,X_1(t),a)
-\int _{-\infty}^{+\infty}\int _0^{t}L _2(s,a)
{\bf \rm d}_{s,a}\nabla _2^-f_v(s,X_1(s),a)\nonumber\\
&&+\int_0^t{\nabla _1^-}{\nabla _2^-}f(s,X_1(s),X_2(s))d<\hskip-4pt M_1,M_2\hskip-4pt>_s.\ \ a.s.
\end{eqnarray}
\end{thm}

Now assume there exists a curve $x_2=b(x_1)$ with left derivative
${{d^-}\over{dx_1}}b(x_1)$ being locally bounded and $b(X_1(t))$
is a semi-martingale. Let $x_2^*=x_2-b(x_1)$ and
$g(t,x_1,x_2^*)=f(t,x_1,x_2^*+b(x_1))$. We can have a generalized
It${\hat {\rm o}}$'s formula in terms of $X_1(s)$ and $X_2^*(s)$
similar to (\ref{f6}). Let $L_2^*(t,a)$ be the local time of
$X_2^*(t)$. In particular, the following result obtained in
\cite{peskir3} can be derived from Theorem \ref{cfr11} as a
special case of our theorem.
\begin{cor}
Assume $f:R^2\to R$ is left continuous and locally bounded and there exists a continuous curve $x_2=b(x_1)$  such that

(i)  the left derivative ${{d^-}\over{dx_1}}b(x_1)$ exists and is locally bounded and
$b(X_1(t))$ is a semi-martingale;

 (ii) $f(x_1,x_2)$ is twice differentiable with continuous second order derivatives ${{\partial^{2}}\over {\partial{x_i}\partial{x_j}}}f$ $(i,j=1,2)$ in regions $x_2\leq b(x_1)$ and $x_2\geq b(x_1)$ respectively.

 Then for any two-dimensional continuous semi-martingale $(X_1(t),X_2(t))$
\begin{eqnarray}\label{f7}
&&
f(X_1(t),X_2(t))-f(X_1(0),X_2(0))\nonumber\\
&=&\sum_{i=1}^2\int _0^t\nabla _i ^-f(X_1(s),X_2(s))dX_i(s)+{1\over 2}\sum_{i=1}^2\int_0^t \Delta_i^-f(X_1(s),X_2(s))d<\hskip-4pt X_i\hskip-4pt>_s\nonumber\\
&&+\int_0^t \Big[\nabla_2^- f(X_1(s),b(X_1(s))+)-\nabla_2^- f(X_1(s),b(X_1(s))-) \Big]d L _2^*(s,0)\nonumber\\
&&+\int_0^t{\nabla _1^-}{\nabla _2^-}f(X_1(s),X_2(s))d<\hskip-4pt M_1,M_2\hskip-4pt>_s.\ \ a.s.
\end{eqnarray}
\end{cor}
{\bf Proof}: Formula (\ref{f7}) can be read from (\ref{f6}) by considering
  \begin{eqnarray*}
f_h(x_1,x_2)&=&f(x_1,x_2)+\int_0^{x_1}(\nabla_2 f(y,b(y)-)-\nabla_2 f(y,b(y)+))(x_2-b(y))^+dy,\\
f_v(x_1,x_2)&=&\int_0^{x_1}(\nabla_2 f(y,b(y)+)-\nabla_2 f(y,b(y)-))(x_2-b(y))^+dy,
 \end{eqnarray*}
 and the integration by parts formula.
To verify conditions of Theorem \ref{cfr11} on $f_v$,  first note
\begin{eqnarray*}
\nabla_1^-f_v(x_1,x_2)&=&(\nabla_2 f(x_1,b(x_1)+)-\nabla_2 f(x_1,b(x_1)-))(x_2-b(x_1))^+,\\
\nabla_2^-f_v(x_1,x_2)&=&\int_0^{x_1}(\nabla_2 f(y,b(y)+)-\nabla_2 f(y,b(y)-))1_{\{x_2> b(y)\}}dy,\\
\nabla_1^-\nabla_2^-f_v(x_1,x_2)&=&(\nabla_2 f(x_1,b(x_1)+)-\nabla_2 f(x_1,b(x_1)-))1_{\{x_2 > b(x_1)\}}.
\end{eqnarray*}
It's trivial to prove that $\nabla_1^-\nabla_2^-f_v(x_1,x_2+b(x_1))$ is of locally bounded variation in $x_2$. To see $\nabla_2^-f_v(x_1,x_2+b(x_1))$ is of locally bounded variation for $x_2$, for any partition $-N=x_2^0<x_2^1<\cdots<x_2^n=N$,
\begin{eqnarray*}
&&\sum\limits_{i=0}^{n-1}|\nabla_2^-f_v(x_1,x_2^{i+1}+b(x_1))-\nabla_2^-f_v(x_1,x_2^{i}+b(x_1))|\\
&\leq&\sum\limits_{i=0}^{n-1}\int_0^{x_1}|\nabla_2f(y,b(y)+)-\nabla_2f(y,b(y)-)|1_{\{x_2^i+b(x_1)\leq b(y)\leq x_2^{i+1}+b(x_1)\}} dy\\
&\leq&\int_0^{x_1}|\nabla_2f(y,b(y)+)-\nabla_2f(y,b(y)-)|1_{\{-N+b(x_1)\leq b(y)\leq N+b(x_1)\}} dy\\
&<&\infty.
\end{eqnarray*}
In order to prove that $\nabla_1^-\nabla_2^-f_v(x_1,x_2+b(x_1))=(\nabla_2 f(x_1,b(x_1)+)-\nabla_2 f(x_1,b(x_1)-))1_{\{x_2 > 0\}}$ and $\nabla_1^-f_v(x_1,x_2+b(x_1))=(\nabla_2 f(x_1,b(x_1)+)-\nabla_2 f(x_1,b(x_1)-))x_2^+$ are of locally bounded variation in $x_1$, we only need to prove that $\nabla_2 f(x_1,b(x_1)+)-\nabla_2 f(x_1,b(x_1)-)$ is of locally bounded variation in $x_1$.
This is true, because for $x_2^*>0$
\begin{eqnarray*}
&&{D^-\over {D{x_1}}} \nabla_2 f(x_1,x_2^*+b(x_1))\\
&=&\nabla_1\nabla_2 f(x_1,x_2^*+b(x_1))+\nabla_2\nabla_2 f(x_1,x_2^*+b(x_1)){d^-\over {dx_1}}b(x_1).
\end{eqnarray*}
So as $x_2^*\to 0+$,
\begin{eqnarray*}
&&{D^-\over {D{x_1}}}  \nabla_2 f(x_1,x_2^*+b(x_1))\\
&\to& \nabla_1\nabla_2 f(x_1, b(x_1)+)+\nabla_2^2f(x_1, b(x_1)+){d^-\over {dx_1}}b(x_1)\\
&=&{d ^-\over {d{x_1}}} \nabla_2 f(x_1,b(x_1)+).
\end{eqnarray*}
It follows that $\nabla_2f(x_1,b(x_1)+)$ is of locally bounded variation. Similarly, one can prove that $\nabla_2f(x_1,b(x_1)-)$ is also of locally bounded variation.
 $\hfill\diamond$\\
 \bigskip

{\bf Acknowledgement} \vskip5pt We would like to acknowledge
partial financial supports to this project by the EPSRC research
grants GR/R69518 and GR/R93582. CF would like to thank the
Loughborough University development fund for its financial
support. It is our great pleasure to thank N. Eisenbaum, D.
Elworthy, Y. Liu, Z. Ma, S. Peng, G. Peskir, A. Truman, J. A. Yan,
M. Yor and W. Zheng for helpful discussions. We would like to
thank G. Peskir and N. Eisenbaum for invitation to the
mini-workshop of local time-space calculus with applications in
Oberwolfach May 2004 where the results of this paper were
announced; to S. Peng for invitation to speak at the 9-th Chinese
mathematics summer school (Weihai) 2004; to F. Gong to the
workshop on stochastic analysis in Chinese Academy of Sciences
2004 and M. Chen to the workshop on stochastic processes and
related topics in Beijing 2004. We would
like to thank the referee for careful reading of the manuscript and pointing out
an error in the early version of the paper and other useful suggestions.

\end{document}